\documentclass[12pt]{amsart}
\usepackage{amsfonts}
\usepackage{amssymb}
\theoremstyle{definition}

\theoremstyle{remark}

\begin{document}

\ \ \ \ \ \ \ \ \ \ \ \ \ \ \ \ \ \ \ \ \ \ \ \ \ \ \ \ \ \ \ \ \ \ \ \ \ \ \ \ \ \ \ \ \ \ \ \ \ \ \ \ \ \ \ \ \ \ \ \ \ \ 
Bull. Math. Soc. Sc. Math. Roumanie

\ \ \ \ \ \ \ \ \ \ \ \ \ \ \ \ \ \ \ \ \ \ \ \ \ \ \ \ \ \ \ \ \ \ \ \ \ \ \ \ \ \ \ \ \ \ \ \ \ \ \ \ \ \ \ 
Tome 38 (86) Nr. 3-4, 1994-1995

\vspace{1.2in}

\centerline{\large\bf Note on Holomorphic }
\centerline{\bf\large Transormations Preserving the }
\centerline{\bf\large Bochner curvature tensor }

\vspace{0.5in}
\centerline{ by }

\centerline{\large Ognian KASSABOV
\footnote{Partially supported by the Ministry of Education of Bulgaria, Contract MM413/94}}

\vspace{0.5in}

\centerline{\bf Abstract }

\vspace{0.2in}
In \cite{B} D. Blair conjectured, that any holomorphic transformation of 
a Kaehler manifold, which preserves the Bochner curvature tensor is a homothety.
He also gives a result in this direction. In this note we give the affirmative
answer to the Blair's conjecture for the case $\dim_RM = 4$.

\vspace{0.6in}
\noindent
{\large\bf 1  Preliminaries}

\vspace{0.2in}
\noindent
Let $M$ be a $2n$-dimensional Kaehler manifold with  metric tensor g and
complex structure $J$. Denote by $R,\, S,\, \tau$ the curvature tensor, the Ricci 
tensor (of type (1,1) as well as of type (0,2)) and the scalar curvature 
of $M$, respectively. Then the Bochner curvature tensor $B$ 
is defined by 
$$
	 B(x,y)z =  R(x,y)z - \frac 1{2(n+2)}\{S(y,z)x - S(x,z)y + g(y,z)Sx - g(x,z)Sy  
$$
$$
	 +S(Jy,z)Jx - S(Jx,z)Jy + g(Jy,z)SJx - g(Jx,z)SJy - 2S(Jx,y)Jz 
$$
$$
	-2g(Jx,y)SJz\}+\frac{\tau}{4(n+1)(n+2)} 
	\{g(y,z)x - g(x,z)y +g(Jy,z)Jx
$$
$$
	-g(Jx,z)Jy - 2g(Jx,y)Jz \} \ .
$$
It is well known, that $B$ has the algebraic properties of the curvature tensor,
see e.g. \cite{TV}. Moreover, if $\{ e_i,\, Je_i \ i=1,\hdots,n \}$ is an
orthonormal basis of $T_pM$ for a point $p \in M$, then
$$
	\sum_{i=1}^n B(e_i,Je_i)x=0       \leqno (1)
$$
holds good for any $x \in T_pM$, see e.g. \cite{TV}.

\vspace{0.3in}
\noindent
{\large\bf 2 Proof of Blair's conjecture for $n=2$}

\vspace{0.3in}
\noindent
{\bf Theorem 2.1.} {\it Let $M$ be a 4-dimensional Kaehler manifold with
nonvanishing Bochner curvature tensor $B$. Then any holomorphic
transformation $f$ of $M$, which preserves $B$ is a homothety.}

\vspace{0.1in}
\noindent
{\bf Proof:} Let \ $x,\, y \in T_pM$, such that \ $B(x,Jy)\ne 0$. We choose an
orthonormal basis $\{ e_1,e_2,Je_1,Je_2 \}$ \ of \ $T_pM$, which diagonalize
the symmetric $J$-invariant endomorphism $B(x,Jy)J$, i.e.
$$
	B(x,Jy)Je_i = \lambda_ie_i   \leqno (2)
$$
for $i=1,2$. Hence, denoting by $B$ the Bochner curvature tensor of type (0,4)
too, we have
$$
	B(x,Jy,Je_1,e_1) + B(x,Jy,Je_2,e_2) = \lambda_1 +\lambda_2
$$
and because of (1) we obtain $ \lambda_1 +\lambda_2=0$. For a vector $x$ and
a tensor field $T$ on $M$ let $\bar x$ and $\overline T$ be defined by $\bar x = f_*x$
and $\overline T = f^*T$, respectively. Since $f$ is holomorphic and preserves the
Bochner tensor, (2) implies
$$
	\overline B(\bar x,J\bar y)J\bar e_i = \lambda_i\bar e_i   \leqno (3)
$$
and hence
$$
	\sum_{i=1}^2\overline B(\bar x,J\bar y,J\bar e_i,\bar e_i) = \sum_{i=1}^2\lambda_i\bar g(\bar e_i,\bar e_i) \ .   \leqno (4)
$$
On the other hand, since $f$ preserves the Bochner tensor, (1) gives
$$
	\overline B(\bar e_1,J\bar e_1)\bar x + \overline B(\bar e_2,J\bar e_2)\bar x = 0   
$$ 
and consequently
$$
	\overline B(\bar e_1,J\bar e_1,J\bar y,\bar x) + \overline B(\bar e_2,J\bar e_2,J\bar y,\bar x) = 0 \ . \leqno (5)   
$$
From (4) and (5) we obtain 
$\lambda_1\bar g(\bar e_1,\bar e_1) + \lambda_2\bar g(\bar e_2,\bar e_2) = 0$.
Since $\lambda_1+\lambda_2 = 0$ and $(\lambda_1,\lambda_2) \ne (0,0)$, it follows
$\bar g(\bar e_1,\bar e_1) = \bar g(\bar e_2,\bar e_2)$. On the other hand, from (3) we have
$$
	\overline B(\bar x,J\bar y,J\bar e_i,\bar e_j) = \lambda_i\bar g(\bar e_i,\bar e_j)   
$$
and hence, using the algebraic properties of $\overline B$ and $\lambda_1 + \lambda_2 = 0$, we obtain:
$$
	\lambda_1\bar g(\bar e_1,\bar e_2)   =
	\overline B(\bar x,J\bar y,J\bar e_1,\bar e_2) = 
	\overline B(\bar x,J\bar y,J\bar e_2,\bar e_1) 
	= \lambda_2\bar g(\bar e_1,\bar e_2) = -\lambda_1\bar g(\bar e_1,\bar e_2)   \ ,
$$
i.e. $\lambda_1\bar g(\bar e_1,\bar e_2) = - \lambda_1\bar g(\bar e_1,\bar e_2) $.
Since $\lambda_1 \ne 0 $ we get $ g(\bar e_1,\bar e_2) = 0$. Consequently 
$\bar g = \mu g$ at $p$, i.e. $f$ is a homothety at $p$. Hence $\mu$ is a constant on $M$,
because as it is easy to see, a conformal, non-homothetic change of a K\"ahler metric 
destroys the K\"ahler property.

\vspace{0.6in}

\ \ \ \ \ \ \ \ \ \ \ \ \ \ \ \ \ \ \ \ \ \ \ \ \ \ \ \ \ \ \ \ \ \ \ \ \ \ \ \ \ \ \ \ \ \ \ \ \ \ \ \ \ \ \ \ \ \ \ \ \ \ 
{\bf Higher Transport School}

\ \ \ \ \ \ \ \ \ \ \ \ \ \ \ \ \ \ \ \ \ \ \ \ \ \ \ \ \ \ \ \ \ \ \ \ \ \ \ \ \ \ \ \ \ \ \ \ \ \ \ \ \ \ \ \ \ \ \ \ \ \ 
{\bf BBTY"T. Kableschkov"}

\ \ \ \ \ \ \ \ \ \ \ \ \ \ \ \ \ \ \ \ \ \ \ \ \ \ \ \ \ \ \ \ \ \ \ \ \ \ \ \ \ \ \ \ \ \ \ \ \ \ \ \ \ \ \ \ \ \ \ \ \ \ 
{\bf Section of Mathematics}

\ \ \ \ \ \ \ \ \ \ \ \ \ \ \ \ \ \ \ \ \ \ \ \ \ \ \ \ \ \ \ \ \ \ \ \ \ \ \ \ \ \ \ \ \ \ \ \ \ \ \ \ \ \ \ \ \ \ \ \ \ \ 
{\bf Slatina, 1574 Sofia}

\ \ \ \ \ \ \ \ \ \ \ \ \ \ \ \ \ \ \ \ \ \ \ \ \ \ \ \ \ \ \ \ \ \ \ \ \ \ \ \ \ \ \ \ \ \ \ \ \ \ \ \ \ \ \ \ \ \ \ \ \ \ 
{\bf BULGARIA}

\vspace{0.2in}
\noindent
Received May 5, 1995

\end{document}